\documentclass[11pt,a4paper]{article}
\usepackage{amsmath, amsfonts, amssymb,}
\usepackage{a4wide}
\usepackage{parskip}
\usepackage{enumitem}
\usepackage{xcolor}
\usepackage{indentfirst}
\usepackage{pstricks}
\usepackage{latexsym}
\usepackage{amsthm}
\usepackage{cite}

\newcommand{\e}{\epsilon}

\newcommand{\al} {\alpha}

\newcommand{\de} {\delta}

\newcommand{\Ga} {\Gamma}
\newcommand{\Om} {\Omega}
\newcommand{\De} {\Delta}
\newcommand{\la} {\lambda}

\newcommand{\noi} {\noindent}

\newcommand{\ds} {\displaystyle}

\newcommand{\RR}{{\mathbb R}}

\newtheorem{thm}{Theorem}[section]
\newtheorem{rem}{Remark}[section]
\newtheorem{lem}{Lemma}[section]
\newtheorem{prop}{Proposition}[section]

\numberwithin{equation}{section}
\def\proof{\noindent{\textbf{Proof. }}}

\begin{document}
	\setlength{\abovedisplayskip}{3pt}
	\setlength{\belowdisplayskip}{3pt}
	\date{}
	{\vspace{0.01in}
		\title{Critical exponent Neumann problem with Hardy-Littlewood-Sobolev nonlinearity }
		
		\author{  Jacques Giacomoni$^{\,1}$ \footnote{e-mail: {\tt jacques.giacomoni@univ-pau.fr}}, \ Sushmita Rawat$^{\,2}$\footnote{e-mail: {\tt sushmita.rawat1994@gmail.com }},  \
			and \  K. Sreenadh$^{\,2}$\footnote{
				e-mail: {\tt sreenadh@maths.iitd.ac.in}} \\
			 $^1\,${\small Universit\'e  de Pau et des Pays de l'Adour, LMAP (UMR E2S-UPPA CNRS 5142) }\\ {\small Bat. IPRA, Avenue de l'Universit\'e F-64013 Pau, France}\\  
			$^2\,${\small Department of Mathematics, Indian Institute of Technology Delhi,}\\
			{\small	Hauz Khaz, New Delhi-110016, India } }


		\maketitle
		\begin{abstract}
		\noi In this article, we study the Brezis–Nirenberg type problem of nonlinear Choquard equation with Neumann boundary condition
		\begin{equation*}
			\begin{aligned}
				-\Delta u &= \lambda \alpha(x) u + \left(\int\limits_{\Omega}\frac{u(y)^{2^*_{\mu}}}{|x-y|^{\mu}}\;dy\right)u^{2^*_\mu-1}, \;\;\text{in} \; \Omega,\\
				\frac{\partial u}{\partial \nu} &= 0\;\; \text{on} \; \partial\Omega,
			\end{aligned}
			\end{equation*}
		where $\Omega$ is a bounded domain in $\mathbb{R}^N  (N\geq 4)$, $\nu$ is the unit outer normal to $\partial \Om$  and $\mu \in (0, N)$. According to the parameter $\lambda$ we prove necessary and sufficient conditions for the existence and non-existence of positive weak solutions to the problem. The proof is based on variational arguments.
		\medskip

		
		\medskip
		

	\end{abstract}
\newpage
\section{Introduction}	
	The purpose of this article is to study a class of Choquard equation with a linear perturbation and Neumann boundary conditions. We consider the following problem:
		\begin{equation*}(P)\; \left\{
			\begin{array}{cc}				
				-\Delta u = \ds\la \al(x) u + \left(\int\limits_{\Om}\frac{u(y)^{2^*_\mu}}{|x-y|^{\mu}}\;dy\right)u^{2^*_\mu-1}, \;\;\text{in} \; \Omega,\\
				\dfrac{\partial u}{\partial \nu} = 0\;\; \text{on} \; \partial\Omega,
			\end{array}\right.
		\end{equation*}	
where $\Omega$ is a bounded domain of $\mathbb{R}^{N}$ having a smooth boundary, $N \geq 4$ with $\lambda \in \RR$, $\al \in C^\infty(\overline{\Om})$ and $\nu$ denoting the unit outer normal to $\partial \Om$. Here $2^{*}_{\mu }= \frac{2N-\mu}{N-2}$ is the critical exponent in the sense of Hardy-Littlewood-Sobolev inequality recalled in Proposition \ref{littlewoodsobolevi}. 

The study of the Neumann boundary conditions with Laplacian operators has been an active area of research for several decades. A considerable body of literature is available for problems involving both sub-critical and critical nonlinearities. Numerous scholars have explored various aspects of this problem, such as the existence, uniqueness, and regularity of solutions, along with their qualitative behavior and applications in physics. One of the earliest works in this area is the paper by Lin, Ni and Takagi \cite{Lin_Ni_Takagi}, where the authors obtained the existence of the least energy solutions for the sub-critical case. Further, in \cite{Ni_Takagi1, Ni_Takagi2} authors proved that the least energy solutions attain only one maximum, and it is achieved at exactly one point on the boundary $\partial\Om$ which maximizes the mean curvature of boundary. These works laid the foundation for subsequent research, which focused on extending the study to more general cases and exploring the behavior of the solutions in greater details. For the critical case, Adimurthi and Yadava \cite{adimurthi_yadava} established existence results. This work is followed by several other papers that investigated the properties of ground state solutions and behavior of solutions under different conditions, see for reference \cite{adi_man1, adi_man2, adi_pac_yad, x.j.wang, z.q.wang1, z.q.wang2}, also  concerning the effect of mixed boundary conditions one can refer to \cite{grossi_pacella} and references therein.

One challenging problem that has not been investigated so far is the study of elliptic equations involving nonlocal terms, as Choquard type equations,  with Neumann boundary conditions. The Choquard equation has attracted a lot of attention in recent years due to its applications in the study of models in quantum mechanics, Bose-Einstein condensation, and nonlinear optics. This type of problem has numerous applications in physical models. In particular, one of the  applications of Choquard equations was given by Pekar\cite{pekar}. He proposed the following problem
\begin{equation}\label{choq}
	-\Delta u + V(x) u = \left(\;\int\limits_{\RR^3}\frac{u^2(y)}{|x-y|}\;dy\right)u \text{ in } \RR^3
\end{equation}
 for the modeling of quantum polaron. Lieb \cite{choqlieb} investigated it in the context of an  approximation to Hartree-Fock theory of one-component plasma. The author proved the existence and uniqueness of \eqref{choq}, up to translations, of the ground state. For a further state of the art on the subject, readers can refer to \cite{Moroz4, M.yang, tuhina_sreenadh, goel_rawat_sreenadh} and references therein.
 
	In the present paper, we will investigate a class of Neumann problems with critical Choquard nonlinear terms and obtain a Brezis-Nirenberg type existence result, that refers to the study of elliptic perturbed equations with critical growth nonlinearities. Precisely, in their pioneering work, Brezis and Nirenberg \cite{brez_nir} investigated nontrivial solutions to a nonlinear elliptic equation of the form:
	\begin{equation*}
			- \De u= \la u + u^{2^*-1} \;\text{ in }\; \Om,
			\end{equation*}
		with Dirichlet boundary conditions $(u = 0 \; \text{ on }\; \partial\Om)$. Here $\Omega\subset \mathbb{R}^N$, $N \geq 4$ is a bounded domain, $2^* = \frac{2N}{N-2}$ denotes the Sobolev critical exponent, $\lambda$ is a constant parameter. They proved the existence of solutions for certain values of $\la > 0$ by exhibiting Palais-Smale sequences whose energy levels are strictly below the first critical level, ensuring compactness. 
		The case of Neumann boundary type conditions was also considered in the context of Brezis-Nirenberg problems.
In particular, assuming  a flatness condition on $\partial \Om$ and for a suitable range of $\la >0$, Adimurthi and Yadava in \cite{adimurthi_yadava}, established existence and nonexistence results for this type of problem.
The Brezis-Nirenberg problem with Hartree type nonlinearities was also investigated. In this regard Gao and Yang in \cite{M.yang} established some existence results for a class of Choquard with Dirichlet boundary conditions. Moreover, in \cite{tuhina_sreenadh}, authors studied the nonlocal counterpart of this problem and obtained various results such as existence, multiplicity, regularity, and nonexistence results. 

In this paper, we aim to contribute to the existing literature by investigating the Brezis-Nirenberg problem in case of Choquard nonlinearity with Neumann boundary conditions that has not been explored before. Specifically, following the approach in \cite{adimurthi_yadava}, we will show sharp conditions for the existence and non-existence of positive solutions to this class of problems. To prove the existence of a minimizer for problem $(P)$ and according to the boundary conditions, we require a suitable version of the inequality of Cherrier\cite[Theorem 1]{cherrier1}, which states that for every $\e > 0$ there exists a constant $C_\e > 0$ such that
\begin{equation*}
	\left( \int\limits_{\Om}| u|^{2^*}\;dx\right)^{\frac{2}{2^*}} \leq  \left( \frac{2^{\frac{2}{n}}}{S} + \e \right)\int\limits_{\Om}|\nabla u|^2\;dx + C(\e)\int\limits_{\Om}|u|^2\;dx,
\end{equation*}	
where $S$ is the best Sobolev constant for the embedding $H^1_0(\Om)$ into $L^{2^{*}}(\Om)$. Following ideas of Aubin\cite{aubin}, we extend this result in the presence of Choquard nonlinearity (see Lemma \ref{Flem3.1}). We highlight that this extension is not straightforward and is of independent interest.
Further, we utilize this inequality to minimize the Sobolev quotient associated with our problem and establish compactness of our Sobolev functional within a range determined by the best constant $S_H$ (defined in \eqref{Feq2.8}) for the Sobolev embedding $H^1_0(\Om)$ into $L^{2^*}(\Om)$.

We now state our main result. Let $\al$ change sign in $\Om$ and $\int\limits_{\Om} \al(x)\,dx <0$. Let $\la(\al)$ be the unique real number such that the problem
\begin{equation}\label{Feq1.2}
	\begin{aligned}
	-\De\phi &= \la(\al)\al(x)\phi \;\;\text{in} \; \Omega,\\
	\frac{\partial \phi}{\partial \nu} &= 0\;\; \text{on} \; \partial\Omega,
\end{aligned}
\end{equation}
admits a positive  solution. Then we prove the following existence result: 
	\begin{thm}\label{Fthm1.1}
		Assume that 
		\begin{enumerate}
			\item[1.] $\al$ changes sign in $\Om$ and $\int\limits_{\Om} \al(x)\,dx <0$. Let $\la(\al)$ be given by \eqref{Feq1.2}.
			\item[2.] There exists a $x_0 \in \partial \Om$ such that $\al(x_0) >0$ and $\partial \Om$ is flat of order $k > \frac{6N-\mu}{2N-\mu}$ at $x_0$.
		\end{enumerate}
Then, problem $(P)$ has a positive solution $u \in C^2(\overline{\Om})$ if and only if $\la \in (0, \la(\al))$.	\end{thm}

Throughout the paper, we make use of the following notations:
\begin{itemize}
	\item $u^{\pm} := \max\{\pm u, 0\}$ 
	\item For any $u \in H^{1}(\mathcal{A})$, we set
	\begin{equation*}
		\|u\|_{0, \mathcal{A}}^{2\cdot2^*_\mu}:= \iint\limits_{\mathcal{A} \times \mathcal{A}} \frac{(u(x))^{2^*_\mu}(u(y))^{2^*_\mu}}{|x-y|^ \mu}\,dxdy.
	\end{equation*}
    \item  $B_R(a)$ stands for the open ball with center $a$ and radius
    $R$ and 
    $B_R^+(a)= B_R(a) \cap \{x_N > 0\}$.
	\item The letters $C$, $C_i$ denote various positive constants possibly different in various places.
\end{itemize} 
The paper is organized as follows: In Section 2, we provide some preliminary results. In Section 3, we prove the main result of the paper by showing the Choquard version of the Cherrier inequality and by establishing accurate estimates on the Sobolev quotient (see Lemma \ref{lem-quotient}).

	\section{Preliminaries}
	We recall some definitions of function spaces and results that will be required in later sections.
	Consider the functional space $H^1(\Om)$ as the usual Sobolev space defined as
	\begin{equation*}
		H^1(\Om) = \left\lbrace u\in L^2(\Om): \int\limits_{\Om} |\nabla u|^2\,dx < \infty \right\rbrace,
	\end{equation*}
	endowed with the norm 
	\begin{equation*}
		\|u\|^2_{H^1(\Om)} = \int\limits_{\Om} |\nabla u|^2\,dx + \int\limits_{\Om} | u|^2\,dx.
	\end{equation*}		
From \cite{adimurthi_yadava}, we recall
\begin{prop}
	Given $\alpha \in C^\infty(\overline{\Omega})$ such that $\alpha$ changes sign in $\Omega$ and $\ds\int_{\Omega} \alpha(x)\;dx <0$, let $\la(\al)$ be the unique real number given by \eqref{Feq1.2}. Then, we have
\begin{enumerate}
	\item[1.] for all $\la \in (0, \la(\al)),$
\begin{equation*}
	\left( \int\limits_{\Om} |\nabla u|^2\,dx - \la\int\limits_{\Om}\al(x)  u^2\,dx\right)^{\frac{1}{2}}
\end{equation*}
defines an equivalent norm on $H^1(\Om)$.
\item[2.] Let $u \in H^2(\Om) \cap C(\overline{\Om})$ be such that
\begin{equation*}
	\begin{array}{cc}
		\la \notin (0, \la(\al))\\
		\Delta u + \la \al u \neq 0,\\
		\Delta  u + \la \al u \leq 0, \dfrac{\partial u}{\partial \nu} = 0,\\
	\end{array}
\end{equation*}
then $u$ can not be positive.
\end{enumerate}
\end{prop}
	\begin{prop}\label{littlewoodsobolevi}
		\textbf{(Hardy-Littlewood-Sobolev inequality)}: Let $t$, $r > 1$ and $0 < \mu < N$ with $\frac{1}{t} + \frac{\mu}{N} + \frac{1}{r} = 2$, $f \in L^t(\mathbb{R}^N)$ and $h \in L^r(\mathbb{R}^N)$. Then there exists a sharp constant $C(t, r, \mu, N)$ independent of $f$, $h$ such that
		\begin{equation*}
			\iint\limits_{\mathbb{R}^{2N}} \dfrac{f(x)h(y)}{|x-y|^ \mu}\,dxdy \leq C(t, r, \mu, N)\|f\|_{L^t(\mathbb{R}^N)}\|h\|_{L^r(\mathbb{R}^N)}.
		\end{equation*}
		\end{prop}
	 In particular, let $f = h = |u|^t$ 
	then by Hardy-Littlewood-Sobolev inequality we see that,
	$$    \|u\|_{0, \RR^{N}}^{2\cdot2^*_{\mu}} :=                          \iint\limits_{\RR^{2N}}\frac{|u(x)^t|u(y)|^t}{|x-y|^{\mu}}dxdy$$
	is well defined if $|u|^t \in L^\nu(\RR^N)$ with $\nu =\frac{2N}{2N-\mu}>1$.
	Sobolev embedding theorems yield
	\begin{equation*}
		\frac{2N-\mu}{N} \leq t \leq \frac{2N-\mu}{N-2}.
	\end{equation*}
	From this, for $u \in H^1(\RR^N) $ one has
	$$    \|u\|_{0, \RR^{N}}^{2} \leq C(N,\mu)^\frac{1}{2_\mu^\ast}  |u|_{2^*}^2  .   $$
	We are looking for nontrivial solutions for problem $(P)$ which are defined as, $u \in H^1(\Om)$ such that
\begin{equation*}
	\int\limits_{\Om}|\nabla u|^2\,dx = \la \int\limits_{\Om}\al(x) u^2\,dx + \iint\limits_{\Omega\times\Omega}\frac{(u(y))^{2^{*}_{\mu}}(u(x))^{2^{*}_{\mu}-1}\phi(x)}{|x-y|^\mu}\,dxdy
\end{equation*}
for any $\phi \in H^1(\Om)$. 
The energy functional associated with the problem $(P)$ is ${J : H^1(\Omega) \rightarrow \mathbb{R}}$\; defined as,
\begin{equation*}
 J(u)= \frac{1}{2} \int\limits_{\Om} |\nabla u|^2\,dx -\frac{\la}{2} \int\limits_{\Omega} \al(x)u(x)^{2} \,dx -\frac{1}{2\cdot2^*_{\mu}}\|u\|_{0, \Om}^{2\cdot2^*_\mu}.
\end{equation*}
From the embedding results, we know that $H^1(\Om)$ is continuously embedded in $L^{p}(\Om)$ when $1 \leq p \leq 2^{*}$. Also the embedding is compact for $1 \leq p < 2^{*}$, but not for the case $p = 2^{*}$.                                  
Accordingly, we denote $S_H$ (independent of $\Omega$) the best constant for the embedding $H^1_0(\Om)$ into $L^{2^{*}}(\Om)$ as
\begin{equation}\label{Feq2.8}
S_H = \inf\limits_{u\in H^1(\Om)\backslash \{0\}}\left\lbrace\int\limits_{\Om} |\nabla u|^2\,dx:  \|u\|_{0,\Om}^{2\cdot2^{*}_{\mu}} = 1\right\rbrace. 
\end{equation}

\begin{lem} \cite[Lemma 1.2]{M.yang}
The constant $S_H$ is achieved by $u$ if and only if $u$ is of the form
$C\left( \frac{t}{t^2 + |x-x_0|^2}\right) ^\frac{N-2}{2}$, $x \in \mathbb{R}^N$, 
for some $x_0 \in \mathbb{R}^N, C \text{and}\; t > 0.$ Moreover,
$S_H = \frac{S}{{C(N, \mu)}^\frac{1}{2_{\mu}^{*}}}$, where $S$ is the best Sobolev constant.
\end{lem}

Consider the family of functions ${U_\epsilon}$, 
defined by
\begin{equation} \label{minmimizer}
U_\epsilon = \epsilon^{-\frac{N-2}{2}}u^{*}\left( \frac{x}{\epsilon}\right), \; x\in \mathbb{R}^{N}, \epsilon > 0,
\end{equation}
\begin{equation*}
u^{*}(x) = \overline{u}\left( \frac{x}{{S}^\frac{1}{2}}\right) , \; \overline{u}(x) = \frac{\tilde{u}(x)}{\|\tilde u\|_{L^{2^{*}}(\mathbb{R}^N)}} \; \text{and}\; \tilde{u}(x) = \alpha(\beta^2 + |x|^2)^{-\frac{N-2s}{2}},
\end{equation*}
with $\alpha > 0$ and $\beta > 0$ are fixed constants.
Then for each $\epsilon > 0, \; U_\epsilon $ satisfies
\begin{equation*}
-\Delta u = |u|^{2^{*}-2}u \quad in \; \mathbb{R}^N,
\end{equation*}
and the equality,
\begin{equation*}
\int\limits_{\mathbb{R}^{N}} |\nabla U_\epsilon(x)|^2\,dx = \int\limits_{\mathbb{R^N}}|U_\epsilon|^{2^{*}} = {S}^\frac{N}{2}. 
\end{equation*}

\textbf{Flatness condition}\\
Let $x_0 \in \partial\Om$, after translation and rotation, we assume that $x_0 =0$ and there exists $R>0$ and a smooth function $\rho: B_R(0) \cap \{x_N =0\} \to \RR$ such that
\begin{equation*}
	\begin{aligned}
		\rho(0) = 0&, \; \nabla\rho(0) = 0\\
		\Om \cap B_R(0) &= \{x \in  B_R(0): x_N > \rho(x')\}\\
		\partial\Om \cap B_R(0) &= \{x \in B_R(0): x_N = \rho(x')\}
	\end{aligned}
\end{equation*}
where $x' = (x_1, x_2, \cdots, x_{N-1}, 0)$.
\begin{rem}
	We say that $\partial\Om$ is of order $k$ at $0$ if $\rho(x') = O(|x'|^k)$, as $|x'|\to 0$.
\end{rem}

\section{Proof of theorem \ref{Fthm1.1}}
In this section we begin by examining the problem with mixed Dirichlet and Neumann boundary conditions. We then derive estimates for this problem. Using this, we deduce a corresponding result for our initial problem.\\
Let $\Om \subset \RR^N$ be a bounded domain with smooth boundary. Let $\Ga_0$ and $\Ga_1$ be disjoint sub-manifolds of $\partial\Om$ such that $\partial \Om = \Ga_0 \cup \Ga_1$ and let
\begin{equation*}
	H^1_{\Gamma_0}(\Omega):= \{u \in H^1(\Om): u \equiv 0 \;\text{on}\; \Ga_0\}.
\end{equation*}
Let $a \in L^\infty(\Om)$, $b \in L^\infty(\Ga_1)$ be such that 
\begin{equation*}
		\|u\|:= \left( \int\limits_{\Om} |\nabla u|^2\,dx - \int\limits_{\Om}a(x)  u^2\,dx + \int\limits_{\Ga_1}b(x)  u^2\,dx\right)^{\frac{1}{2}}
\end{equation*}
defines an equivalent norm on $H^1_{\Ga_0}(\Omega))$. For $u \in H^1_{\Ga_0}(\Omega)$, define
\begin{equation*}
	\begin{aligned}
		J(u)&= \frac{1}{2} \|u\|^2 -\frac{1}{2\cdot2^*_{\mu}}\|u\|_{0,\Om}^{2\cdot2^{*}_{\mu}}\\
		Q(u) &= \frac{\|u\|^2}{\|u\|_{0,\Om}^{2}}\\
		S_H(\Ga_0, a, b) &= \inf\limits_{u \in H^1_{\Ga_0}(\Omega)\backslash\{0\}}Q(u).
	\end{aligned}
\end{equation*}
Further, we extend the results of Cherrier \cite{cherrier1, cherrier2} to the Hardy-Littlewood-Sobolev critical case, and obtain the following crucial lemma:
\begin{lem}\label{Flem3.1}
	For every $\e>0$, there exists $C_\e >0$ such that for all $u \in H^1(\Om)$,
\begin{equation*}
	\|u\|_{0,\Om}^2 \leq  \left( \frac{2^{\frac{2^*_\mu-2}{2^*_{\mu}}}}{S_H} + \e \right)\int\limits_{\Om}|\nabla u|^2\;dx + C(\e)\int\limits_{\Om}|u|^2\;dx.
\end{equation*}	
\end{lem}
\proof First, we prove the result for the half-space $\RR^N_{+}$. Let $v \in H^1(\RR^N_+)$, defined as
\begin{equation*}
	u(x',x_N) = \begin{cases}
		v(x', x_N) & \;\text{if}\; x_N >0;\\
		v(x', -x_N) & \;\text{if}\; x_N <0,
	\end{cases}
\end{equation*}
where $x= (x',x_N)$. By direct calculation it follows that
\begin{equation}\label{Feq23.2}
	\int\limits_{\RR^N}|\nabla u|^2\;dx = 2 \int\limits_{\RR^N_+}|\nabla v|^2\;dx.
\end{equation}
Similarly, we see that
\begin{equation}\label{Feq23.3}
\|u\|_{0, \RR^{N}}^{2\cdot2^*_{\mu}} = 4\|v\|_{0,\RR^{N}_+}^{2\cdot2^*_{\mu}}.
\end{equation}
Employing Hardy-Littlewood-Sobolev inequality, Sobolev embedding theorem, \eqref{Feq23.2} and \eqref{Feq23.3} we obtain that
\begin{equation}\label{Feqb3.4}
	\begin{aligned}
		\|v\|_{0,\RR^N_+}^2 = \frac{\|u\|_{0,\RR^N}^2}{2^{\frac{2}{2^*_{\mu}}}}
		\leq \frac{1}{2^{\frac{2}{2^*_{\mu}}}S_H} \int\limits_{\,\RR^N}|\nabla u|^{2}dx
		= \frac{2^{\frac{2^*_{\mu}-2}{2^*_{\mu}}}}{S_H}\int\limits_{\,\RR^N_+}|\nabla v|^{2}dx.
	\end{aligned}
\end{equation}
From \eqref{Feqb3.4}, we deduce that for any open set $E$ of $\RR^N_+$ and for any $v\in H^1(\RR^N_+)$ with compact support in $E$ we have that
	\begin{equation}\label{Feqb3.5}
		\|v\|_{0,E}^2\leq \frac{2^{\frac{2^*_\mu-2}{2^*_\mu}}}{S_H}\int_E|\nabla v|^2dx.
	\end{equation}
For the bounded domain case, we  are inspired by the approach of Aubin\cite[Theorem 1.31, Lemma 2.26]{aubin}, where they obtained the result by constructing local charts with suitable properties (via exponential mapping there) for a compact manifold without boundary (\cite[Lemma 2.24]{aubin}). Furthermore, we will adapt some ideas in the proof of \cite[Theorem 2.30]{aubin} concerning manifolds with boundaries to the bounded domain $\Omega$ with smooth boundary case. Precisely given the connectedness and compactness of $\bar{\Omega}$, we infer 
that for any $\epsilon>0$, there exists $\delta=\delta(\epsilon)>0$ and a finite covering $\displaystyle(\Omega_i)_{1\leq i\leq m}$ with $\Omega_i\supset B_{\de}(p_i)$ and $p_i\in \Omega$, or $\Omega_i\supset B^+_{\de}(p_i)$ and $p_i\in \partial\Omega$ and corresponding local charts $(\Omega_i,\phi_i)_{1\leq i\leq m}$, $\phi_i:\, \Omega_i\mapsto  B_i$,  either $B_i=B_{r_i}(0) $ or $B^+_{r_i}(0)$ and verifying :
\begin{equation}\label{Feqb3.6}
	(1+\epsilon')^{\frac{-2^*_\mu}{\mu}}\leq ||D\phi_i(x)||\leq {(1+\epsilon')^{\frac{2^*_\mu}{\mu}}}
\end{equation}
with $\epsilon'= 2^{\frac{2-2^*_{\mu}}{2^*_\mu}}\epsilon\,S_H$.
For that, let us detail the case $p_i\in \partial\Omega$ (for $p_i\in \Omega$ it is even easier, see \cite[Lemma 2.24]{aubin}). 
Up to a rotation of coordinate axis and translation, we can assume that $p_i=0_{\RR^N}$, $x_N=\rho(x')$ is a local parametrization near  $p_i$ of the boundary $\partial\Omega$ and $x_N=0$ is the tangent hyperspace at $p_i$ of $\partial\Omega$ (that implies $\frac{\partial \rho}{\partial x_i}(0)=0$ for $i=1,..., N-1$).
Consider the mapping $\phi\,:\, (x',x_N)\mapsto (x',x_N+\rho(x'))$ that 
satisfies $D\phi(0)=I$. For $r>0$ small enough, $\phi_{|_{B_r^+(0)}}$ is then a diffeomorphism and we can take 
$\phi_i=\phi^{-1}_{|_{\phi(B_r^+(0))}}$ and $\Omega_i=\phi(B_r^+(0))$ and up to reducing $r_i=r$ together with the compactness of $\Omega$, we get \eqref{Feqb3.6}.\\
 Consider a partition of unity $\{ h_i\}_{1\leq i\leq m}\subset C^\infty$ subordinate to the atlas $(\Omega_i, \phi_i)_{1\leq i\leq m}$. Then, 
using \eqref{Feqb3.5},\eqref{Feqb3.6}, one has
\begin{equation*}
	\begin{aligned}
\|u\|_{0,\Omega}&=\|\displaystyle\sum_{i=1}^mh_i u\|_{0,\Omega}\leq 
	\displaystyle\sum_{i=1}^m\left(\;\iint\limits_{\phi_i(\Omega_i)\times\phi_i(\Omega_i)}
\frac{|h_i u(\phi_i^{-1}(x))|^{2^*_\mu}|h_i u(\phi_i^{-1}(y))|^{2^*_\mu}}{|\phi_i^{-1}(x)-\phi_i^{-1}(y)|^\mu}dxdy\right)^{\frac{1}{2\cdot2^*_{\mu}}}\\
&\leq
	(1+\epsilon')^{1/2}
	\displaystyle\sum_{i=1}^m\left(\;\iint\limits_{\phi_i(\Omega_i)\times\phi_i(\Omega_i)}
		\frac{|h_i u(\phi_i^{-1}(x))|^{2^*_\mu}|h_i u(\phi_i^{-1}(y))|^{2^*_\mu}}{|x-y|^\mu}dxdy\right)^{\frac{1}{2\cdot2^*_{\mu}}}\\
		&\leq \frac{2^{\frac{2^*_{\mu}-2}{2\cdot2^*_{\mu}}}}{S_H^{1/2}}
			(1+\epsilon')^{1/2}\displaystyle\sum_{i=1}^m\left(\int_{\phi_i(\Omega_i)}|\nabla h_i u(\phi^{-1}_i(x))|^2dx\right)^{1/2}\\
&\leq \frac{2^{\frac{2^*_{\mu}-2}{2\cdot2^*_{\mu}}}}{S_H^{1/2}}	(1+\epsilon')^{1/2}	\left[\displaystyle\sum_{i=1}^m\left(\int_{\phi_i(\Omega_i)}h_i|\nabla u(\phi_i^{-1}(x))|^2dx\right)^{1/2}\right.\\
&\quad+C \left.	\displaystyle\sum_{i=1}^m\left(\int_{\phi_i(\Omega_i)}u^2(\phi_i^{-1}(x))dx\right)^{1/2}\right].
\end{aligned}
\end{equation*}
This together with the inequality $(A+B)^2\leq (1+\eta)A^2 +C(\eta)B^2$ for any $\eta>0$,
imply that
\begin{equation*}
	\left(\;\iint \limits_{\Omega\times\Omega}\frac{(u(x))^{2^*_\mu}(u(y))^{2^*_\mu}}{|x-y|^ \mu}\,dxdy \right)^{\frac{1}{2^*_\mu}}\leq
		\left(\frac{2^{\frac{2^*_\mu-2}{2^*_{\mu}}}}{S_H}+\epsilon\right)\int_\Omega|\nabla u|^2dx+C(\epsilon)\int_{\Omega}u^2dx.
\end{equation*}
This completes the proof of Lemma \ref{Flem3.1}. \qed

Using Lemma \ref{Flem3.1}, we obtain the following result.
\begin{lem}
	The following holds:
	\begin{enumerate}
		\item[1.] $S_H(\Ga_0, a, b) >0$;
		\item[2.] Assume $S_H(\Ga_0, a, b) < \frac{S_H}{2^{\frac{2^*_\mu-2}{2^*_{\mu}}}}$, then there exists a $v \geq 0$ such that $S_H(\Ga_0, a, b)= Q(v)$. 
		Further if we define $u_0 = \left( S_H(\Ga_0, a, b)\right) ^\frac{1}{2\cdot2^*_{\mu}-2}v$, then $u_0$ satisfies
		\begin{equation}\label{Feq3.2}
			\begin{aligned}
				-\De u_0 &= a(x)u_0 + \left(\int\limits_{\Om}\frac{u_0(y)^{2^*_\mu}}{|x-y|^{\mu}}\;dy\right)u_0^{2^*_\mu-1} \;\;\text{in}\; \Om,\\
				u_0 &>0,\\
				u_0 &= 0 \;\;\text{on}\; \Ga_0, \;\frac{\partial u_0}{\partial \nu} + bu_0 = 0 \;\;\text{on}\; \Ga_1,
			\end{aligned}
		\end{equation} 
		and $J(u_0) < \ds \left(\frac{1}{2} - \frac{1}{2\cdot2^*_{\mu}} \right)\frac{1}{2^\frac{2^*_{\mu}-2}{2^*_{\mu}-1}}S_H^\frac{2^*_{\mu}}{2^*_{\mu}-1}$. 

	\end{enumerate}
\end{lem}
\proof   By the Hardy-Littlewood-Sobolev inequality and the Sobolev embedding theorem, for all $u \in H^1_{\Ga_0}(\Omega)$, we have that
\begin{equation*}
	\|u\|_{0,\Om}^{2} \leq \frac{\|u\|^2}{S_H},
\end{equation*}
and the proof of 1 follows by the definition of $S_H(\Ga_0, a,b)$.\\
Proof of 2: Consider a minimizing sequence $\{u_n \}$ for $S_H(\Ga_0, a,b)$ such that $\|u\|_{0,\Om}^{2\cdot2^*_{\mu}}=1$. Let for a subsequence, $u_n \rightharpoonup v$ weakly in $H^1_{\Ga_0}(\Omega)$ as $n \to \infty$. \\
Claim 1: $v \equiv 0$. \\
If suppose $v \equiv 0$, that is, $u_n \rightharpoonup 0$, then by the Rellich lemma and Lemma \ref{Flem3.1}, we obtain
\begin{equation*}
	\begin{aligned}
		\lim\limits_{n \to \infty}\int\limits_{\Om}|\nabla u_n|^2\,dx 
		&= S_H(\Ga_0, a, b)\\
		&\leq S_H(\Ga_0, a, b)\left( \frac{2^{\frac{2^*_\mu-2}{2^*_{\mu}}}}{S_H} + \e \right)\lim\limits_{n \to \infty}\int\limits_{\Om}|\nabla u_n|^2\;dx,
	\end{aligned}
\end{equation*}
for every $\e > 0$. Hence, we have, 
\begin{equation*}
	1 \leq S_H(\Ga_0, a, b)\left( \frac{2^{\frac{2^*_\mu-2}{2^*_{\mu}}}}{S_H} + \e \right).
\end{equation*}
This contradicts $S_H(\Ga_0, a, b) < \frac{S_H}{2^{\frac{2^*_\mu-2}{2^*_{\mu}}}}$. Therefore, $v \neq 0$. \\
Claim 2: $Q(v) = S_H(\Ga_0, a, b)$. \\
Let $v_n := u_n -v$ and $v_n  \rightharpoonup 0$ in $H^1_{\Ga_0}(\Omega)$. By Brezis-Lieb lemma \cite[Theorem 1]{BrezLieb}, we have as $n \to \infty$ 
\begin{equation*}
	\|u_n\|^2 = \|v\|^2 + \int\limits_{\Om}|\nabla v_n|^2\;dx + o(1),
\end{equation*}
that is,
\begin{equation}\label{Feq3.3}
	S_H(\Ga_0, a, b) = \|v\|^2 + \int\limits_{\Om}|\nabla v_n|^2\;dx + o(1).
\end{equation}
From \cite[Lemma 2.2]{M.yang} and Lemma \ref{Flem3.1}, we deduce that as $n \to \infty$,
\begin{equation*}
	1 = \|u_n\|_{0,\Om}^{2\cdot2^*_{\mu}} \leq \|v\|_{0,\Om}^{2\cdot2^*_{\mu}} + \left( \frac{2^{\frac{2^*_\mu-2}{2^*_{\mu}}}}{S_H} + \e \right)\int\limits_{\Om}|\nabla v_n|^2\;dx + o(1),
\end{equation*}
for any $\e >0$, that is,
\begin{equation}\label{Feq3.4}
		S_H(\Ga_0, a, b) \leq S_H(\Ga_0, a, b)\|v\|_{0,\Om}^{2\cdot2^*_{\mu}} + \int\limits_{\Om}|\nabla v_n|^2\;dx + o(1).
\end{equation}
From \eqref{Feq3.3} and \eqref{Feq3.4}, we obtain the desired claim. As $Q(v)= Q(|v|)$, we can assume that $v \geq 0$. Also by defining $u_0 = \left( S_H(\Ga_0, a, b)\right)^\frac{1}{2\cdot2^*_{\mu}-2}v$, we see that $u_0$ satisfies problem \eqref{Feq3.2}. Moreover,
\begin{equation*}
\begin{aligned}
		J(u_0) &= \frac{\left( S_H(\Ga_0, a, b)\right) ^\frac{1}{2^*_{\mu}-1}}{2}\|v\|^2- \frac{\left( S_H(\Ga_0, a, b)\right) ^\frac{2^*_{\mu}}{2^*_{\mu}-1}}{2\cdot2^*_{\mu}}\|v\|_{0,\Om}^{2\cdot2^*_{\mu}}
	& < \left(\frac{1}{2} - \frac{1}{2\cdot2^*_{\mu}} \right)\frac{S_H^\frac{2^*_{\mu}}{2^*_{\mu}-1}}{2^\frac{2^*_{\mu}-2}{2^*_{\mu}-1}}.
\end{aligned}
\end{equation*}

\begin{lem}\label{lem-quotient}
	Assuming $\al\in C(\overline{\Om})$ and there exists a point $x_0\in \partial\Om$ where $\al(x_0) > 0$ and the boundary is flat of order $k > \frac{6N-\mu}{2N-\mu}$ at $x_0$. Then for every $\la > 0$, we have
	\begin{equation*}
S_H(\la\al) < \frac{S_H}{2^{\frac{2^*_\mu-2}{2^*_{\mu}}}},
	\end{equation*} 
	where $S_H(\la\al)= S_H(\phi, \la\al, 0)$. 
\end{lem}
\proof 
Assuming no loss of generality, we can set $x_0 = 0$, so that $0\in \partial\Om$, $\al(0)> 0$, and $\partial \Om$ is flat at $0$ of order $k > \frac{6N-\mu}{2N-\mu}$. Let $\rho: B_R(0)\cap \{x:x_N = 0\} \to \RR$ be the function used in the definition of flatness. For any $u \in H^1(\Om)$, we have
\begin{equation*}
	Q(u) = \frac{\ds\int\limits_{\Om} |\nabla u|^2\,dx - \la\int\limits_{\Om}\al(x)  u^2\,dx}{\|u\|_{0,\Om}^{2}}.
\end{equation*}
Let $\phi\in C_c^{\infty}(B_{R/2}(0))$ be such that $\phi$ is radial and $\phi = 1$ on $B_{R/4}(0)$ and for each $\epsilon > 0$, let  $u_\epsilon$ be defined as 
\begin{equation*}
	u_\e(x) = \phi(x)U_\epsilon(x) =\frac{\phi(x)}{\left( \e+|x|^2\right)^{\frac{N-2}{2}} },
\end{equation*}
where $U_\epsilon$ is defined in \eqref{minmimizer}.
We claim that as $\e \to 0$,
\begin{equation}\label{Feq3.5}
	Q(u_\e) < 2^{\frac{2-2^*_\mu}{2^*_{\mu}}} S_H,
\end{equation}
which implies the lemma. The proof of \eqref{Feq3.5} follows in several steps. Without loss of generality, we can assume that $\rho \geq 0$. For non-positive $\rho$ the estimate \eqref{Feq3.5} follows exactly as in the case of positive $\rho$ with slight modifications. Set
\begin{equation*}
	\Sigma := \{x \in B_{R/2}(0) : 0 < x_N < \rho(x')\}.
\end{equation*}
From \cite{adimurthi_yadava}, we have the following estimate
\begin{equation}\label{Feq3.6}
\begin{aligned}
		\int\limits_{\Om} |\nabla u_\e|^2\,dx &= \frac{1}{2}\int\limits_{B_{R}(0)} |\nabla u_\e|^2\,dx - \int\limits_{\Sigma} |\nabla u_\e|^2\,dx\\
	&= \frac{\|\nabla u\|_{L^2}^2}{2\e^{\frac{N-2}{2}}}\left( 1+ O(\e^{\frac{N-2}{2}}) + O(\e^{\frac{k-1}{2}})\right),
\end{aligned}
\end{equation}
and 
\begin{equation}\label{Feq3.7}
	\begin{aligned}
		\int\limits_{\Om} | u_\e|^2\,dx &= \frac{1}{2}\int\limits_{B_{R}(0)} |u_\e|^2\,dx - \int\limits_{\Sigma} |u_\e|^2\,dx\\
		&= \begin{cases}
			\dfrac{\|u\|_{L^2}^2}{2\e^{\frac{N-4}{2}}}\left( 1+ O(\e^{\frac{N-4}{2}}) + O(\e^{\frac{k-1}{2}})\right) & \;\text{if}\; N \geq 5 \vspace{0.3cm}\\ 
			\dfrac{\omega}{2}|\log \e|+ O(1) + O(\e^{\frac{k-1}{2}}) & \;\text{if}\; N = 4,\\
		\end{cases}
	\end{aligned}
\end{equation}
where $\omega$ is the area of $S^3$ and $\frac{\|\nabla u\|_{L^2}^2}{\|u\|_0^2} = S_H$. Next, we estimate the Choquard term
\begin{equation}\label{Feq3.8}
	\|u_\e\|_{0,\Om}^{2\cdot2^*_{\mu}} = \frac{1}{4}\iint\limits_{B_{R}(0) \times B_{R}(0)} \frac{(u_\e)^{2^*_\mu}(u_\e)^{2^*_\mu}}{|x-y|^ \mu}\,dxdy - \iint\limits_{\Sigma \times \Sigma} \frac{(u_\e)^{2^*_\mu}(u_\e)^{2^*_\mu}}{|x-y|^ \mu}\,dxdy -2 \iint\limits_{\Om \times \Sigma} \frac{(u_\e)^{2^*_\mu}(u_\e)^{2^*_\mu}}{|x-y|^ \mu}\,dxdy.
\end{equation}
Consider
\begin{equation}\label{Feq3.9}
	\begin{aligned}
		\|u_\e\|_{0,B_{R}(0)}^{2\cdot2^*_{\mu}} &\geq \iint\limits_{B_{\frac{R}{4}}(0) \times B_{\frac{R}{4}}(0)}\frac{dxdy}{\left( \e + |x|^2\right)^{\frac{2N-\mu}{2}} \left( \e + |y|^2\right)^{\frac{2N-\mu}{2}}|x-y|^{\mu} }\\
      & = \frac{1}{\e^{\frac{2N-\mu}{2}}}\iint\limits_{B_{\frac{R}{4\sqrt{\e}}}(0) \times B_{\frac{R}{4\sqrt{\e}}}(0)}\frac{dxdy}{\left( 1 + |x|^2\right)^{\frac{2N-\mu}{2}} \left( 1 + |y|^2\right)^{\frac{2N-\mu}{2}}|x-y|^{\mu} }	\\
      & = \frac{1}{\e^{\frac{2N-\mu}{2}}}\left\lbrace \|u\|_{0,\RR^N}^{2\cdot2^*_{\mu}} -2 \mathbb{D}- \mathbb{E}  \right\rbrace ,
	\end{aligned}
\end{equation}
where 
\begin{equation*}
	\mathbb{D} = \iint\limits_{\RR^N \backslash B_{\frac{R}{4\sqrt{\e}}}(0) \times B_{\frac{R}{4\sqrt{\e}}}(0)}\frac{dxdy}{\left( 1 + |x|^2\right)^{\frac{2N-\mu}{2}} \left( 1 + |y|^2\right)^{\frac{2N-\mu}{2}}|x-y|^{\mu} }
\end{equation*}
and
\begin{equation*}
	\mathbb{E} = \iint\limits_{\RR^N \backslash B_{\frac{R}{4\sqrt{\e}}}(0) \times \RR^N \backslash B_{\frac{R}{4\sqrt{\e}}}(0)}\frac{dxdy}{\left( 1 + |x|^2\right)^{\frac{2N-\mu}{2}} \left( 1 + |y|^2\right)^{\frac{2N-\mu}{2}}|x-y|^{\mu} }.
\end{equation*}
By utilizing both the Hardy-Littlewood-Sobolev inequality and the beta function, we have
\begin{equation}\label{Feq3.10}
	\begin{aligned}
	\mathbb{D} &\leq C(N,\mu)\left[ \int\limits_{B_{\frac{R}{4\sqrt{\e}}}(0)}\frac{dx}{\left( 1 + |x|^2\right)^{N}}\right] ^{\frac{2N-\mu}{2N}}
	\left[ \int\limits_{\RR^N \backslash B_{\frac{R}{4\sqrt{\e}}}(0)}\frac{dy}{\left( 1 + |y|^2\right)^{N}}\right] ^{\frac{2N-\mu}{2N}}\\
	&\leq C(N,\mu)\left[\int\limits_0^{\infty}\frac{r^{N-1}\,dr}{(1+r^2)^N} - \int\limits_{\frac{R}{4\sqrt{\e}}}^{\infty}\frac{r^{N-1}\,dr}{(1+r^2)^N} \right]^{\frac{2N-\mu}{2N}}\left[ \int\limits_{\frac{R}{4\sqrt{\e}}}^{\infty}\frac{r^{N-1}\,dr}{(1+r^2)^N} \right]^{\frac{2N-\mu}{2N}}\\
		&\leq C(N,\mu) \left[\frac{\Ga(N/2)\Ga(N/2)}{\Ga(N)} \right]^{\frac{2N-\mu}{2N}}\e^{\frac{2N-\mu}{4}}, 
	\end{aligned}
\end{equation}
and 
\begin{equation}\label{Feq3.11}
	\begin{aligned}
		\mathbb{E} \leq C	\left[ \int\limits_{\RR^N \backslash B_{\frac{R}{4\sqrt{\e}}}(0)}\frac{dy}{\left( 1 + |y|^2\right)^{N}}\right] ^{\frac{2N-\mu}{N}} 
		\leq C \e^{\frac{2N-\mu}{2}}. 
	\end{aligned}
\end{equation}
Putting together \eqref{Feq3.10} and \eqref{Feq3.11} in \eqref{Feq3.9}, we get
\begin{equation}\label{Feq3.12}
	\|u_\e\|_{0,B_{R}(0)}^{2\cdot2^*_{\mu}} \geq \frac{\|u\|_{0,\RR^N}^{2\cdot2^*_{\mu}}}{\e^{\frac{2N-\mu}{2}}} \left( 1 - O(\e^{\frac{2N-\mu}{4}}) \right). 
\end{equation}
Next, with the help of Hardy-Littlewood-Sobolev inequality and \cite{adimurthi_yadava}, we infer that
\begin{equation}\label{Feq3.13}
\begin{aligned}
		\iint\limits_{\Sigma \times \Sigma} \frac{(u_\e)^{2^*_\mu}(u_\e)^{2^*_\mu}}{|x-y|^ \mu}\,dxdy & \leq 
		C(N,\mu)\left[\int\limits_{\frac{R}{4} < |x| < \frac{R}{2}}\frac{(\phi(x))^{2^*}-1}{(\e+ |x|^2)^N}\;dx + \int\limits_{\Sigma }\frac{dx}{(\e+ |x|^2)^N}\right]^{\frac{2N-\mu}{N}}\\
		& = \left(O(1) + O(\e^{\frac{k-1-N}{2}})\right)^{\frac{2N-\mu}{N}}\\
		& = \begin{cases}
			\dfrac{1}{\e^{\left(\frac{N+1-k}{2} \right) \left(\frac{2N-\mu}{N} \right) }} \left(O(1) + O(\e^{\frac{N+1-k}{2}})\right) & \;\text{if}\; k \leq N+1 \vspace{0.25cm}\\
			O(1) + O(\e^{\frac{k-1-N}{2}}) & \;\text{if}\; k >N+1, 
		\end{cases}
\end{aligned}
\end{equation}
and
\begin{equation}\label{Feq3.14}
	\begin{aligned}
		\iint\limits_{\Om \times \Sigma} \frac{(u_\e)^{2^*_\mu}(u_\e)^{2^*_\mu}}{|x-y|^ \mu}\,dxdy & \leq C(N,\mu)\left[ \int\limits_{\Om } (u_\e(x))^{2^*}\,dx\right]^{\frac{2N-\mu}{2N}}\left[ \int\limits_{\Sigma } (u_\e(x))^{2^*}\,dx\right]^{\frac{2N-\mu}{2N}}\\
		& = \left[ \frac{k_2^{\frac{2^*}{2}}}{2\e^{\frac{N}{2}}}\left(1 + O(\e^{\frac{N}{2}}) + O(\e^{\frac{k-1}{2}}) \right)\right]^ {\frac{2N-\mu}{2N}} \left(O(1) + O(\e^{\frac{k-1-N}{2}})\right)^{\frac{2N-\mu}{2N}}\\
		& =\begin{cases} \ds \frac{k_2^{\frac{2^*_\mu}{2}}\left(O(1) +  O(\e^{\frac{k-1}{2}})+ O(\e^{\frac{N+1-k}{2}}) \right)}{2^{\frac{2^*_{\mu}}{2^*}}\e^{\frac{2N-\mu}{4}\left[1 +\frac{N+1-k}{N} \right]} } & \;\text{if}\; k \leq N+1, \vspace{0.25cm}\\
			\ds\frac{k_2^{\frac{2^*_\mu}{2}}}{2^{\frac{2^*_{\mu}}{2^*}}\e^{\frac{2N-\mu}{4}}}\left( O(1) + O(\e^{\frac{k-1-N}{2}}) + O(\e^{\frac{N}{2}})\right) & \;\text{if}\; k >N+1.
		\end{cases}
	\end{aligned}
\end{equation}
Thus depending upon the range of $k$, we consider the following cases:\\
\textbf{Case I:} $k \leq N+1$\\
It follows from \eqref{Feq3.12}, \eqref{Feq3.13}, \eqref{Feq3.14} and \eqref{Feq3.8}
\begin{equation}\label{Feq3.15}
	\begin{aligned}
		\|u_\e\|_{0,\Om}^{2\cdot2^*_{\mu}} &\geq \frac{\|u\|_{0,\RR^N}^{2\cdot2^*_{\mu}}}{4\e^{\frac{2N-\mu}{2}}} \left( 1 - O(\e^{\frac{2N-\mu}{4}}) \right) - \frac{1}{\e^{\left(N+1-k \right) \left(\frac{2N-\mu}{2N} \right) }} \left(O(1) + O(\e^{\frac{N+1-k}{2}})\right)\\ & \qquad-\frac{2k_2^{\frac{2^*_\mu}{2}}\left(O(1) + O(\e^{\frac{k-1}{2}})+ O(\e^{\frac{N+1-k}{2}}) \right)}{2^{\frac{2^*_{\mu}}{2^*}}\e^{\frac{2N-\mu}{4}\left[1 +\frac{N+1-k}{N} \right] }}\\
		& = \frac{\|u\|_{0,\RR^N}^{2\cdot2^*_{\mu}}}{4\e^{\frac{2N-\mu}{2}}} \left( 1 - O(\e^{\frac{2N-\mu}{4}}) - O(\e^{(k-1)\frac{2N-\mu}{4N}})\right).
	\end{aligned}
\end{equation}
Now choose $R > 0$ and $\al_0 > 0$ such that $\al(x) \geq \al_0$ for all $x \in B_{R}(0) \cap \overline{\Om}$, then 
\begin{equation*}
	Q(u_\e)  \leq \frac{\ds\int\limits_{\Om} |\nabla u_\e|^2\,dx - \la\al_0\int\limits_{\Om} u_{\e}^2\,dx}{\|u_\e\|_{0,\Om}^{2}}.
\end{equation*}
From \eqref{Feq3.6}, \eqref{Feq3.7} and \eqref{Feq3.15}, we get
\begin{equation}\label{Feq3.16}
	\begin{aligned}
		Q(u_\e) &\leq \frac{S_H}{2^{\frac{2^*_{\mu}-2}{2^*_{\mu}}}}\frac{\left( 1+ O(\e^{\frac{N-2}{2}}) + O(\e^{\frac{k-1}{2}})) \right) }{ \left( 1 - O(\e^{\frac{2N-\mu}{4}}) - O(\e^{(k-1)\frac{2N-\mu}{4N}})\right)} \\ 
		& - \frac{C\la\al_0\, \e}{2^{\frac{2^*_{\mu}-2}{2^*_{\mu}}}\left( 1 - O(\e^{\frac{2N-\mu}{4}}) - O(\e^{(k-1)\frac{2N-\mu}{4N}})\right)}\times \begin{cases}
		\left( 	1+ O(\e^{\frac{N-4}{2}}) + O(\e^{\frac{k-1}{2}}))\right)   &\;\text{if}\;  N \geq 5 \vspace{0.3cm}\\ 
			\left( |\log \e|+ O(1) + O(\e^{\frac{k-1}{2}})\right) &\;\text{if}\;  N = 4\\
		\end{cases}
		\\
		& = \frac{S_H}{2^{\frac{2^*_{\mu}-2}{2^*_{\mu}}}} + O(\e^{(k-1)\frac{2N-\mu}{4N}}) + O(\e^{\frac{N-2}{2}}) - \begin{cases}
			C\la \e & \;\text{if}\; N \geq 5,\\
			C\la \e|\log(\e)| & \;\text{if}\; N =4.
		\end{cases}
	\end{aligned}
\end{equation}
Hence, from \eqref{Feq3.16}, we infer that $Q(u_\e) < \frac{S_H}{2^{\frac{2^*_{\mu}-2}{2^*_{\mu}}}}$ whenever $	k > \frac{6N-\mu}{2N-\mu} $.\\
\textbf{Case II:} $k > N+1$\\
We proceed similar to Case I, from \eqref{Feq3.12}, \eqref{Feq3.13}, \eqref{Feq3.14} and \eqref{Feq3.8}, we have
\begin{equation}\label{Feq3.17}
	\begin{aligned}
		\|u_\e\|_{0,\Om}^{2\cdot2^*_{\mu}} &\geq \frac{\|u\|_{0,\RR^N}^{2\cdot2^*_{\mu}}}{4\e^{\frac{2N-\mu}{2}}} \left( 1 - O(\e^{\frac{2N-\mu}{4}}) \right) - \left(O(1) + O(\e^{\frac{k-1-N}{2}})\right)\\ & \;\qquad-\frac{2k_2^{\frac{2^*_\mu}{2}}\left(O(1) + O(\e^{\frac{N}{2}}) + O(\e^{\frac{k-1-N}{2}}) \right)}{2^{\frac{2^*_{\mu}}{2^*}}\e^{\frac{2N-\mu}{4}}}\\
		& = \frac{\|u\|_{0,\RR^N}^{2\cdot2^*_{\mu}}}{4\e^{\frac{2N-\mu}{2}}} \left( 1 - O(\e^{\frac{2N-\mu}{4}})\right).
	\end{aligned}
\end{equation}
Further, by \eqref{Feq3.6}, \eqref{Feq3.7} and \eqref{Feq3.17}, we get for some constant $C>0$,
\begin{equation}\label{Feq3.18}
	\begin{aligned}
		Q(u_\e) &\leq \frac{S_H}{2^{\frac{2^*_{\mu}-2}{2^*_{\mu}}}}\frac{\left( 1+ O(\e^{\frac{N-2}{2}}) \right) }{ \left( 1 - O(\e^{\frac{2N-\mu}{4}}) \right)} \\
		& - \quad\frac{C\la\al_0\, \e}{2^{\frac{2^*_{\mu}-2}{2^*_{\mu}}}\left( 1 - O(\e^{\frac{2N-\mu}{4}}) \right)}\times \begin{cases}
			 \left( 1+ O(\e^{\frac{N-4}{2}}) \right)  &\;\text{if}\;  N \geq 5 \vspace{0.3cm}\\ 
			\left( |\log \e|+ O(1) + O(\e^{\frac{k-1}{2}})\right)  &\;\text{if}\;  N = 4
		\end{cases}\\
		& = \frac{S_H}{2^{\frac{2^*_{\mu}-2}{2^*_{\mu}}}} + O(\e^{\frac{2N-\mu}{4}}) + O(\e^{\frac{N-2}{2}}) - \begin{cases}
			C\la \e & \;\text{if}\; N \geq 5,\\
			C\la \e|\log(\e)| & \;\text{if}\; N =4.
		\end{cases}
	\end{aligned}
\end{equation}
Hence, from \eqref{Feq3.18}, we infer that $Q(u_\e) < \frac{S_H}{2^{\frac{2^*_{\mu}-2}{2^*_{\mu}}}}$.\qed

		\end{document}